\def\titlep{Recursive boson system in the Cuntz algebra ${\cal O}_{\infty}$}
\documentclass[11pt]{article}
\usepackage{graphicx,ifthen}
\usepackage{amssymb}

 at12pt

\newcommand{\sdag}{\scriptsize \dag}

\newcommand{\qed}{\hbox{\rule[-2pt]{3pt}{6pt}}}
\newcommand{\qedh}{\hfill\qed \\}

\newcommand{\vv}{\vspace{.3in}}

\newtheorem{Thm}{Theorem}[section]

\newtheorem{rem}[Thm]{Remark}

\newtheorem{ex}[Thm]{Example}
\newtheorem{defi}[Thm]{Definition}
\newtheorem{lem}[Thm]{Lemma}

\newtheorem{prop}[Thm]{Proposition}

\newtheorem{fact}[Thm]{Fact}

\newcommand{\ww}{\vv\noindent}

\def\cal#1{\mathcal #1}
\def\con{{\cal O}_{N}}
\def\coni{{\cal O}_{\infty}}
\def\edot{=1,\ldots,N}
\def\pr{{\it Proof.}\quad}

\def\scm#1{S({\bf C}^{N})^{\otimes #1}}

\def\co#1{{\cal O}_{#1}}
\def\ltn{l_{2}({\bf N})}

\def\disp#1{{\displaystyle #1}}
\def\brl{branching law}
\def\bfsnl{{\rm BFS}_{N}(\Lambda)}

\addtocounter{footnote}{1}
\def\sftt#1{
\setcounter{equation}{0}
\addtocounter{footnote}{1}
\section{#1}
}

\def\ssft#1{\subsection{#1}}

\begin{document}
%
%
\def\autherp{Katsunori Kawamura}
\def\emailp{e-mail: kawamura@kurims.kyoto-u.ac.jp.}
\def\addressp{College of Science and Engineering Ritsumeikan University,\\
1-1-1 Noji Higashi, Kusatsu, Shiga 525-8577,Japan
}

\def\infw{\Lambda^{\frac{\infty}{2}}V}
\def\zhalfs{{\bf Z}+\frac{1}{2}}
\def\ems{\emptyset}
\def\pmvac{|{\rm vac}\!\!>\!\! _{\pm}}
\def\vac{|{\rm vac}\rangle _{+}}
\def\dvac{|{\rm vac}\rangle _{-}}
\def\ovac{|0\rangle}
\def\tovac{|\tilde{0}\rangle}
\def\expt#1{\langle #1\rangle}
\def\zph{{\bf Z}_{+/2}}
\def\zmh{{\bf Z}_{-/2}}
\def\brl{branching law}
\def\bfsnl{{\rm BFS}_{N}(\Lambda)}
\def\scm#1{S({\bf C}^{N})^{\otimes #1}}
\def\mqb{\{(M_{i},q_{i},B_{i})\}_{i=1}^{N}}
\def\zhalf{\mbox{${\bf Z}+\frac{1}{2}$}}
\def\zmha{\mbox{${\bf Z}_{\leq 0}-\frac{1}{2}$}}
\newcommand{\mline}{\noindent
\thicklines
\setlength{\unitlength}{.1mm}
\begin{picture}(1000,5)
\put(0,0){\line(1,0){1250}}
\end{picture}
\par
 }
\def\sd#1{#1^{\sdag}}
\def\dlim{{\cal D}\mbox{-}\lim}
\def\dsum{{\cal D}\mbox{-}\sum}
\def\fsum{F\mbox{-}\sum}
\def\bx{\mbox{\boldmath$x$}}

%
%
%
\setcounter{section}{0}
\setcounter{footnote}{0}
\setcounter{page}{1}
\pagestyle{plain}

%
%
\title{\titlep}
\author{\autherp\thanks{\emailp}
\\ 
\addressp
}
\date{}
\maketitle

%
%
\begin{abstract}
Bosons and fermions are often written by elements of other algebras.
M. Abe gave a recursive realization of the boson by
formal infinite sums of the canonical generators of 
the Cuntz algebra ${\cal O}_{\infty}$. 
We show that such formal infinite sum always makes sense on a certain dense 
subspace of any permutative representation of ${\cal O}_{\infty}$.
In this meaning, we can regard as if 
the algebra ${\cal B}$ of bosons was a unital 
$*$-subalgebra of ${\cal O}_{\infty}$ on a given permutative representation
by keeping their unboundedness.
By this relation, we compute branching laws arising from restrictions
of representations of ${\cal O}_{\infty}$ on ${\cal B}$.
For example, it is shown that the Fock representation of ${\cal B}$ is given as 
the restriction of the standard representation of ${\cal O}_{\infty}$
on ${\cal B}$.
\end{abstract}


%
%
\sftt{Introduction}
\label{section:first}

Bosons and fermions are not only important in physics but also 
interesting in mathematics.
Studies of their algebras spurred the development of 
the theory of operator algebras \cite{BR}.
Representations of bosons are used to describe representations of 
several algebras \cite{AOS, KQS, MZ}.
Bosons and fermions are often written by elements of other algebras
and such descriptions are useful for several computations.
For example, the boson-fermion correspondence \cite{MJD,Oko01} is well-known.
It is shown that bosons and fermions are corresponded
as operators on the infinite wedge representation of fermions.

%
%
\ssft{Motivation}
\label{subsection:firstone}
In our previous paper \cite{AK1}, 
we have presented a recursive construction of the CAR
(=canonical anticommutation relation) algebra for fermions in terms of 
the Cuntz algebra $\co{2}$ and 
shown that it may provide us a useful tool to study properties 
of fermion systems by using explicit expressions in terms of generators 
of the algebra. 
Let $s_{1},s_{2}$ be the canonical generators of $\co{2}$, that is,
they satisfy that
\[s_{i}^{*}s_{j}=\delta_{ij}I\quad(i,j=1,2),\quad s_{1}s_{1}^{*}+
s_{2}s_{2}^{*}=I.\]
Let $\zeta$ be the linear map on $\co{2}$ defined by
$\zeta(x)\equiv s_{1}xs_{1}^{*}-s_{2}xs_{2}^{*}$ for $x\in \co{2}$.
We recursively define the family $\{a_{1},a_{2},a_{3},\ldots\}$ by
\[a_{1}\equiv s_{1}s_{2}^{*},\quad a_{n}\equiv \zeta(a_{n-1})\quad(n\geq 2).\]
Then $\{a_{n}:n\in {\bf N}\}$ satisfies that
\[a_{n}a_{m}^{*}+a_{m}^{*}a_{n}=\delta_{nm}I,\quad
a_{n}a_{m}+a_{m}a_{n}=a_{n}^{*}a_{m}^{*}+a_{m}^{*}a_{n}^{*}=0
\quad(n,m\in {\bf N})\]
where ${\bf N}=\{1,2,3,\ldots\}$.
We call such $\{a_{n}:n\in {\bf N}\}$ by a {\it recursive fermion system (=RFS)} 
in $\co{2}$.
From this description, the C$^{*}$-algebra ${\cal A}$ generated by
fermions is embedded into $\co{2}$ as a C$^{*}$-subalgebra with common unit:
\[{\cal A}\equiv C^{*}\langle\{a_{n}:n\in {\bf N}\}\rangle
\hookrightarrow \co{2}\]
Furthermore ${\cal A}$ coincides with the fixed-point subalgebra of $\co{2}$
with respect to the $U(1)$-gauge action.
Because every $a_{n}$ is written as a polynomial in
the canonical generators of $\co{2}$ and their $*$-conjugates,
their description is very simple and it is easy to compute
the restriction $\pi|_{{\cal A}}$ of a representation $\pi$ 
of $\co{2}$ on ${\cal A}$.
By using the RFS,
we obtain several new results about fermions \cite{AK06,AK4,AK05,AK02R}.
For example,
assume that $({\cal H},\pi)$ is a $*$-representation of $\co{2}$
with a cyclic vector $\Omega$.
If $\Omega$ satisfies
$\pi(s_{1})\Omega=\Omega$, then $\pi|_{{\cal A}}$ is 
equivalent to the Fock representation of ${\cal A}$ with the vacuum $\Omega$.
If $\Omega$ satisfies $\pi(s_{1}s_{2})\Omega=\Omega$,
then $\pi|_{{\cal A}}$ is equivalent to the direct sum of 
the infinite wedge representation and the dual infinite wedge
representation of ${\cal A}$ \cite{IWF01}.
In this way, well-known results of fermions are explicitly reformulated by
the representation theory of $\co{2}$.

From this, we speculate that the boson can be also simply written by the 
generators of a certain Cuntz algebra like the RFS,
where the boson means a family $\{a_{n}:n\in {\bf N}\}$ satisfying that
%
%
\begin{equation}
\label{eqn:boson}
a_{n}a_{m}^{*}-a_{m}^{*}a_{n}=\delta_{nm}I,\quad a_{n}a_{m}-a_{m}a_{n}
=a_{n}^{*}a_{m}^{*}-a_{m}^{*}a_{n}^{*}=0
\end{equation}
for each $n,m\in {\bf N}$.
However, 
the boson 
is always represented as a family of unbounded operators on a Hilbert space.
Hence the $*$-algebra generated by $\{a_{n}:n\in {\bf N}\}$
never be a $*$-subalgebra of any C$^{*}$-algebra.
On the other hand, the C$^{*}$-algebra approach of boson is well-known
as the CCR algebra (CCR = canonical commutation relations, 
see $\S$ 5.2 in \cite{BR}).
Because the CCR algebra is not a separable C$^{*}$-algebra,
it is impossible to embed it into any Cuntz algebra
as a C$^{*}$-subalgebra.
From these problems, it seems that
a RFS-like description of bosons by any Cuntz algebra is impossible. 

%
%
\ssft{Recursive boson system}
\label{subsection:firsttwo}
In spite of such problems, Mitsuo Abe gave a ``formal" realization
of the boson by the canonical generators of the Cuntz algebra $\coni$ 
in 2006 as follows.
Let $\{s_{n}:n\in {\bf N}\}$ be the canonical generators of $\coni$,
that is,
\[s_{i}^{*}s_{j}=\delta_{ij}I\quad(i,j\in {\bf N}),\quad
\sum_{i=1}^{k}s_{i}s_{i}^{*}\leq I\quad(\mbox{for any }k\in {\bf N}).\]
Define the family $\{a_{n}:n\in {\bf N}\}$ of formal sums by
%
%
\begin{equation}
\label{eqn:rbsthree}
a_{1}\equiv \disp{\sum_{m=1}^{\infty}\sqrt{m}\,s_{m}s_{m+1}^{*}},
\quad a_{n}\equiv \rho(a_{n-1})\quad(n\geq 2)
\end{equation}
where $\rho$ is the formal canonical endomorphism of $\co{\infty}$ defined by
\[\rho(x)\equiv \sum_{n=1}^{\infty}s_{n}xs_{n}^{*}\quad(x\in\coni).\]
By formal computation, we can verify that $a_{n}$'s 
satisfy (\ref{eqn:boson}) where
we assume that infinite sums can be freely exchanged. 
However infinite sums in these equations {\it do not} 
converge in $\coni$ in general. 
Hence (\ref{eqn:rbsthree}) does not make sense as elements of $\coni$.

We show that the Abe's formal description (\ref{eqn:rbsthree}) can be justified 
as unbounded operators defined on
a certain dense subspace of any permutative representation of $\coni$.
Define $X_{N}\equiv \{1,\ldots,N\}$ for $2\leq N<\infty$ and
$X_{\infty}\equiv {\bf N}$.
Let $\{s_{n}:n\in X_{N}\}$ be the set of canonical generators of $\con$
for $2\leq N\leq \infty$.

%
%
\begin{defi}
\label{defi:firstone} \cite{BJ,DaPi2,DaPi3} 
A representation $({\cal H},\pi)$ of $\con$ is permutative 
if there exists a complete orthonormal basis $\{e_{n}\}_{n\in\Lambda}$
of ${\cal H}$ and a family $f=\{f_{i}\}_{i=1}^{N}$
of maps on $\Lambda$ such that $\pi(s_{i})e_{n}=e_{f_{i}(n)}$
for each $n\in\Lambda$ and $i\edot$.
We call $\{e_{n}\}_{n\in\Lambda}$
and the linear hull ${\cal D}$ of $\{e_{n}\}_{n\in\Lambda}$
by the reference basis and the reference subspace 
of $({\cal H},\pi)$, respectively.
\end{defi}
Remark that for any permutation representation
$({\cal H},\pi)$ of $\con$ with the reference subspace ${\cal D}$,
$\pi(s_{n}){\cal D}\subset {\cal D}$
and  $\pi(s_{n}^{*}){\cal D}\subset {\cal D}$ for each $n$,
but $\pi(x){\cal D}\not\subset {\cal D}$ for $x\in \coni$
in general.

From Definition \ref{defi:firstone} and 
(\ref{eqn:rbsthree}), we immediately obtain the following fact.
%
%
\begin{fact}
\label{fact:mainrbs}
For any permutative representation $({\cal H},\pi)$ of $\coni$,
define the family $\{A_{n}:n\in {\bf N}\}$ of operators on 
the reference subspace ${\cal D}$ of $({\cal H},\pi)$ by
\[
\begin{array}{rl}
A_{1}v\equiv &\sum_{m=1}^{\infty}\sqrt{m}\pi(s_{m}s_{m+1}^{*})v,\\
\\
A_{n}v\equiv &
\sum_{m_{n-1}=1}^{\infty}
\cdots
\sum_{m_{1}=1}^{\infty}
\sum_{m=1}^{\infty}
\sqrt{m}\pi(s_{m_{n-1}}\cdots s_{m_{1}}s_{m}s_{m+1}^{*}
s_{m_{1}}^{*}\cdots s_{m_{n-1}}^{*})v
\end{array}
\]
for $v\in {\cal D}$ and $n\geq 2$.
Then the family $\{A_{n}:n\in {\bf N}\}$ satisfies (\ref{eqn:boson})
on ${\cal D}$.
\end{fact}

\noindent
Infinite sums in Fact \ref{fact:mainrbs} 
are actually finite for each $v\in {\cal D}$.
By comparing Fact \ref{fact:mainrbs} and (\ref{eqn:rbsthree}), 
we see that (\ref{eqn:rbsthree}) is well-defined on
the reference subspace of any permutative representation of $\coni$.
Furthermore, the mapping
%
%
\begin{equation}
\label{eqn:rbsone}
 a_{n}\mapsto A_{n}\quad (n\in {\bf N})
\end{equation}
defines a unital $*$-representation $\pi_{{\cal B}}$ 
of the algebra ${\cal B}$ of bosons on ${\cal D}$, that is,
$\pi_{{\cal B}}(a_{n})\equiv A_{n}$ for each $n$.
In consequence,
we obtain the operation
\[({\cal H},\pi)\mapsto ({\cal D},\pi_{{\cal B}})\]
for any permutative representation $({\cal H},\pi)$
of $\coni$ to the representation $({\cal D},\pi_{{\cal B}})$
of ${\cal B}$.
We call $({\cal D},\pi_{B})$ the {\it restriction of} $({\cal H},\pi)$ 
on ${\cal B}$ and often write it by $({\cal H},\pi|_{{\cal B}})$
for convenience in this paper.
Strictly speaking, this is not a restriction
because ${\cal B}$ is neither a subalgebra of $\coni$ 
nor $\pi(\coni){\cal D}\subset {\cal D}$.
%
%
\begin{rem}
{\rm
If a C$^{*}$-algebra ${\cal A}$ irreducibly acts on a Hilbert space ${\cal H}$,
then any (unbounded) operator on ${\cal H}$
can be written by the strong operator limit of elements of ${\cal A}$ 
on ${\cal H}$.
However such description always depends on the choice of 
representation.
Fact \ref{fact:mainrbs} claims that
the description (\ref{eqn:rbsthree}) 
{\it always} hold on {\it any} permutative representations of $\coni$
nevertheless
there exist infinitely many inequivalent permutative representations of $\coni$
and they are not always irreducible.
}
\end{rem}

%
%
\begin{defi}
The family $\{A_{n}:n\in {\bf N}\}$ 
in Fact \ref{fact:mainrbs} is called the recursive
boson system (=RBS) in $\co{\infty}$ with respect to
a permutative representation $({\cal H},\pi)$ of $\coni$.
\end{defi}

\noindent
We identify $A_{n}$ in Fact \ref{fact:mainrbs} with $a_{n}$. 

Remark that ${\cal B}$ is neither a subalgebra of $\coni$
nor that of the double commutations $\pi(\coni)^{''}$ of $\pi(\coni)$.
However for any permutative representation $({\cal H},\pi)$ of $\coni$,
we obtain a representation of the boson as ${\cal B}$ by the RBS.
In this sense, it seems that ${\cal B}$ is a subalgebra of $\coni$
in special situation: 
\[{\cal B}={\rm Alg}\langle \{a_{n},a_{n}^{*}:n\in{\bf N}\}\rangle
\quad \risingdotseq\quad
\mbox{subalgebra of }\coni.
\]
%
%
\ssft{Representations of bosons arising from
permutative representations of $\coni$}
\label{subsection:firstthree}
We show the significance of the RBS in the representation theory
of operator algebras.
The algebra ${\cal B}$ of bosons always appears with a representation
in theoretical physics.
Especially, the Fock representation plays the most important
role among representations of ${\cal B}$.
It has both the mathematical simple structure and the physical meaning.
By the RBS, we can understand the Fock representation from a viewpoint of 
the representation theory of $\coni$.

First, we explain the notion of branching law.
For a group $G$, if there exists an embedding of $G$ into 
some other group $G^{'}$,
then any representation $\pi$ of $G^{'}$ induces 
the restriction $\pi|_{G}$ of $\pi$ on $G$.
The representation $\pi|_{G}$ is not irreducible in general
even if $\pi$ is irreducible.
If $\pi|_{G}$ is decomposed into the direct sum of
a family $\{\pi_{\lambda}:\lambda\in\Lambda\}$ of irreducible representations 
of $G$, then
the equation
\[\pi|_{G}=\bigoplus_{\lambda\in\Lambda}\pi_{\lambda}\]
is called the {\it branching law} of $\pi$.
The branching law can be also considered for a pair of subalgebra and algebra.
Thanks to the RBS, we can consider (an analogy of) branching laws
of permutative representations of $\coni$ which are restricted on ${\cal B}$. 
%
%
\begin{Thm}
\label{Thm:restriction}
\begin{enumerate}
\item
For $j\geq 1$,
let $({\cal H},\pi_{j})$ be a representation
of $\coni$ with a cyclic vector $\Omega$ satisfying
\[\pi_{j}(s_{j})\Omega=\Omega.\]
Then there exists a dense subspace ${\cal D}_{j}$ of ${\cal H}$
and an action $\eta_{i}$ of ${\cal B}$ on ${\cal D}_{j}$ such that
$\eta_{j}({\cal B})\Omega={\cal D}_{j}$ and 
%
%
\begin{equation}
\label{eqn:etaone}
\eta_{j}(a_{n}a_{n}^{*})\Omega=j\Omega\quad(n\geq 1).
\end{equation}
In particular,
$\eta_{1}$ is the Fock representation of ${\cal B}$ with the vacuum $\Omega$. 
\item
Let $({\cal H},\pi_{12})$ be a representation
of $\coni$ with a cyclic vector $\Omega$ satisfying
\[\pi_{12}(s_{1}s_{2})\Omega=\Omega.\]
Then there exist two subspaces $V_{1}$ and $V_{2}$ of ${\cal H}$
and two actions $\eta_{12}$ and $\eta_{21}$ of ${\cal B}$ on $V_{1}$
and $V_{2}$, respectively
such that $V_{1}\oplus V_{2}$ is dense in ${\cal H}$,
$V_{1}=\eta_{12}({\cal B})\Omega$, $V_{2}=\eta_{21}({\cal B})\Omega^{'}$
for $\Omega^{'}\equiv \pi(s_{2})\Omega$ and the following holds:
\[\left\{
\begin{array}{ll}
\eta_{12}(a_{2n-1})\Omega=\eta_{21}(a_{2n})\Omega^{'}=0,\\
\\
\eta_{12}(a_{2n}^{*}a_{2n})\Omega
=\Omega,\\
\\
\eta_{21}(a_{2n-1}^{*}a_{2n-1})\Omega^{'}=\Omega^{'}
\end{array}
\right.
\quad(n\geq 1).\]
\item
Any two of representations in 
$\{\eta_{j},\eta_{12},\eta_{21}: j\geq 1\}$ of ${\cal B}$
are not unitarily equivalent.
\item
All of representations 
$\{\eta_{j},\eta_{12},\eta_{21}: j\geq 1\}$ of ${\cal B}$ 
are irreducible.
\end{enumerate}
\end{Thm}

\noindent
Every representations of $\coni$ in 
Theorem \ref{Thm:restriction} (i) and (ii) are irreducible
permutative representations.
Hence Theorem \ref{Thm:restriction} shows branching laws of 
representations of $\coni$ restricted on ${\cal B}$:
\[\pi_{j}|_{{\cal B}}=\eta_{j}\quad(j\geq 1),\quad
\pi_{12}|_{{\cal B}}=\eta_{12}\oplus \eta_{21}.\]
By comparison to the fermion case in $\S$ \ref{subsection:firstone},
this result shows that the RBS is very similar to the RFS
in a sense of the representation theory of operator algebras.
This result shows the naturality of 
the description in (\ref{eqn:rbsthree}).

In $\S$ \ref{section:second},
we show permutative representations of $\coni$ and
several representations of ${\cal B}$.
In $\S$ \ref{subsection:secondthree},
we prove Theorem \ref{Thm:restriction}.
In $\S$ \ref{section:third}, we show examples.
In $\S$ \ref{subsection:thirdtwo},
we give an interpretation of representations of bosons in 
Theorem \ref{Thm:restriction} by formal infinite product 
of operators.
%
%
\sftt{Representations and their relations}
\label{section:second}
In order to show Theorem \ref{Thm:restriction},
we introduce several representations of $\coni$ and ${\cal B}$.
After this preparation, we show their relations as the
proof of Theorem \ref{Thm:restriction}. 
%
%
\ssft{Permutative representation of Cuntz algebras}
\label{subsection:secondone}
For $N=2,3,\ldots,+\infty$, 
let $\con$ be the {\it Cuntz algebra} \cite{C}, that is, a C$^{*}$-algebra 
which is universally generated by $s_{1},\ldots,s_{N}$ satisfying
$s_{i}^{*}s_{j}=\delta_{ij}I$ for $i,j\edot$ and
\[\sum_{i=1}^{N}s_{i}s_{i}^{*}=I\quad(\mbox{if } N<+\infty),\quad
\sum_{i=1}^{k}s_{i}s_{i}^{*}\leq I,\quad k= 1,2,\ldots,\quad
(\mbox{if }N = +\infty)\]
where $I$ is the unit of $\con$.
Because $\con$ is simple, that is, there is no
nontrivial closed two-sided ideal,
any homomorphism from $\con$ to a C$^{*}$-algebra is injective.
If $t_{1},\ldots,t_{n}$ are elements of a unital C$^{*}$-algebra
${\cal A}$ such that
$t_{1},\ldots,t_{n}$ satisfy the relations of canonical generators of $\con$,
then the correspondence $s_{i}\mapsto t_{i}$ for $i\edot$
is uniquely extended to a $*$-embedding
of $\con$ into ${\cal A}$ from the uniqueness of $\con$.
Therefore we call such a correspondence 
among generators by an embedding of $\con$ into ${\cal A}$.

Define $X_{N}\equiv \{1,\ldots,N\}$ for $2\leq N<\infty$
and $X_{\infty}\equiv {\bf N}$.
For $N=2,\ldots,\infty$ and $k=1,\ldots,\infty$, 
define the product set $X_{N}^{k}\equiv (X_{N})^{k}$ of $X_{N}$.
Let $\{s_{n}:n\in X_{N}\}$ be the set of canonical generators of $\con$
for $2\leq N\leq \infty$.
%
%
\begin{defi}
\label{defi:first}
For $J=(j_{l})_{l=1}^{k}\in X_{N}^{k}$ with $1\leq k < \infty$,
we write $P_{N}(J)$ the class of representations $({\cal H}, \pi)$ of $\con$ 
with a cyclic unit vector $\Omega\in {\cal H}$
such that $\pi(s_{J})\Omega=\Omega$
and $\{\pi(s_{j_{l}}\cdots s_{j_{k}})\Omega\}_{l=1}^{k}$
is an orthonormal family in ${\cal H}$
where $s_{J}\equiv s_{j_{1}}\cdots s_{j_{k}}$.
Here, $\{\pi(s_{j_{l}}\cdots s_{j_{k}})\Omega\}_{l=1}^{k}$ is 
called the \textit{cycle} of $P_{N}(J)$.
\end{defi}

\noindent
We call the vector $\Omega$ in Definition \ref{defi:first}
by the {\it GP vector} of $({\cal H},\pi)$.
A representation $({\cal H},\pi)$ of $\con$ is called a {\it cycle} 
if there exists $J\in X_{N}^{k}$ for $1\leq k < \infty$ 
such that $({\cal H},\pi)$ belongs to $P_{N}(J)$.
Any permutative representation is uniquely decomposed into
cyclic permutative representations up to unitary equivalence.
For any $J$, $P_{N}(J)$ contains only one unitary equivalence class
\cite{BJ,DaPi2,DaPi3,K1}.
We show properties of $P_{\infty}(j)$ $(j\geq 1)$ 
and $P_{\infty}(12)$ more closely as follows.
%
%
\begin{lem}
\label{lem:permutative}
Let ${\cal T}\equiv \{P_{\infty}(j),P_{\infty}(12):j\geq 1\}$.
\begin{enumerate}
\item
For each $X\in {\cal T}$,
any two representations belonging to $X$ are unitarily equivalent.
\item
Any two of representations in ${\cal T}$ are not unitarily equivalent.
\item
All of representations in ${\cal T}$ are irreducible.
\end{enumerate}
\end{lem}

\noindent
The proof of Lemma \ref{lem:permutative} are given
in Appendix \ref{section:appone}.
From Lemma \ref{lem:permutative} (i),
we use symbols $P_{\infty}(j),P_{\infty}(12)$ as their representatives.

For $2\leq N<\infty$,
let $t_{1},\ldots,t_{N}$ be the canonical generators of $\con$.
Define the representation $(\ltn,\pi)$ of $\con$ by
\[\pi(t_{i})e_{n}\equiv e_{N(n-1)+i}\quad(i=1,\ldots,N,\,n\in {\bf N}).\]
Then $(\ltn,\pi)$ is $P_{N}(1)$ of $\con$.
If we identify $\coni$ with a C$^{*}$-subalgebra
of $\con$ by the embedding of $\co{\infty}$ into $\con$ defined by
%
%
\begin{equation}
\label{eqn:embedding}
s_{(N-1)(k-1)+i}\equiv t_{N}^{k-1}t_{i}\quad(k\geq 1,\,i=1,\ldots,N-1),
\end{equation}
then $(\ltn,\pi|_{\co{\infty}})$ is $P_{\infty}(1)$ of $\co{\infty}$.

%
%
\ssft{Representations of bosons}
\label{subsection:secondtwo}
We summarize several representations of bosons and their properties.
We write ${\cal B}$ the $*$-algebra generated by $\{a_{n}:n\in {\bf N}\}$
which satisfies (\ref{eqn:boson}).
A {\it representation} of ${\cal B}$
is a pair $({\cal H},\pi)$ such that ${\cal H}$ is a complex Hilbert space 
with a dense subspace ${\cal D}$ and 
$\pi$ is a $*$-homomorphism from ${\cal B}$ to 
the $*$-algebra $\{x\in {\rm End}_{{\bf C}}({\cal D}): x^{*}{\cal D}\subset
{\cal D}\}$.
A {\it cyclic vector} of $({\cal H},\pi)$ is a vector $\Omega\in {\cal D}$
such that $\pi({\cal B})\Omega={\cal D}$.
%
%
\begin{defi}
\label{defi:bosonrep}
\begin{enumerate}
\item
For $j\geq 1$,
we write $F_{j}$ the class of representations $({\cal H},\pi)$ of ${\cal B}$
with a cyclic vector $\Omega$ satisfying
$\pi(a_{n}a_{n}^{*})\Omega=j\Omega$ for each $n\in {\bf N}$.
\item
We write $F_{12}$ the class of representations $({\cal H},\pi)$ of ${\cal B}$
with a cyclic vector $\Omega$ satisfying
\[\pi(a_{2n-1})\Omega=0,\quad \pi(a_{2n}^{*}a_{2n})\Omega=\Omega\]
for each $n\in {\bf N}$.
\item
We write $F_{21}$ the class of representations $({\cal H},\pi)$ of ${\cal B}$
with a cyclic vector $\Omega$ satisfying
\[\pi(a_{2n})\Omega=0,\quad \pi(a_{2n-1}^{*}a_{2n-1})\Omega=\Omega\]
for each $n\in {\bf N}$.
\end{enumerate}
\end{defi}
A representation $({\cal H},\pi)$ of ${\cal B}$ is called {\it irreducible}
if the commutant of $\pi({\cal B})$ in ${\cal B}({\cal H})$
is the scalar multiples of $I$.

%
%
\begin{lem}
\label{lem:bosonresult}
Let ${\cal S}\equiv \{F_{j},F_{12},F_{21}:j\geq 1\}$.
\begin{enumerate}
\item
For each $X\in {\cal S}$,
any two representations belonging to $X$ are unitary equivalent.
From this, we can identify a representation belonging to $X\in {\cal S}$ 
with $X$.
\item
Any two of representations in ${\cal S}$
are not unitarily equivalent.
\item
All of representations in ${\cal S}$ are irreducible.
\end{enumerate}
\end{lem}
Lemma \ref{lem:bosonresult} is proved in Appendix \ref{section:apptwo}.
We consider the case $j=1$ in Definition \ref{defi:bosonrep} (i).
Then 
$\pi(a_{n}a_{n}^{*})\Omega=\Omega$ for each $n$.
From this, $\pi(a_{n}^{*}a_{n})\Omega=0$.
This implies that $\pi(a_{n})\Omega=0$ for each $n$.
Because $\Omega$ is a cyclic vector,
$F_{1}$ is the Fock representation of ${\cal B}$ with the vacuum $\Omega$.

In this study, we became the first to find 
$F_{j},F_{12},F_{21}$ from the computation of 
branching laws of permutative representations of $\coni$.
After finding the equations of bosons and the vector $\Omega$,
we found the conditions of $F_{j},F_{12},F_{21}$
without using permutative representations of $\coni$.

%
%
\ssft{Proof of Theorem \ref{Thm:restriction}}
\label{subsection:secondthree}
Before the proof, we summarize basic relations of the RBS
$\{a_{n}:n\in {\bf N}\}$ and the canonical generators
$\{s_{n}:n\in {\bf N}\}$ of $\coni$.
From (\ref{eqn:rbsthree}), the following holds
on the reference subspace of any permutative representation
of $\coni$:
\[s_{m}a_{n}=a_{n+1}s_{m},\quad s_{m}a_{n}^{*}=a_{n+1}^{*}s_{m}\quad
(n,m\in {\bf N}),\]
\[\rho(x)s_{i}=s_{i}x\quad(x\in \coni,\,i\in{\bf N}).\]

\noindent
(i)
Fix $j\geq 1$.
First,
we see that $({\cal H},\pi_{j})$ is $P_{\infty}(j)$ with the GP vector $\Omega$.
We simply write $\pi_{j}(s_{n})$ by $s_{n}$ for each $n$.
Define
\[{\cal D}_{j}\equiv {\rm Lin}\langle\{s_{J}\Omega:J\in {\bf N}^{*}\}\rangle\]
where ${\bf N}^{*}\equiv \coprod_{l\geq 1}{\bf N}^{l}$.
Then ${\cal D}_{j}$ is the reference subspace.
We simply write $\{a_{n}:n\in {\bf N}\}$ the RBS on $P_{\infty}(j)$
and ${\cal B}$ the algebra generated by them.
From (\ref{eqn:rbsthree}), 
\[a_{n}a_{n}^{*}=\sum_{K\in{\bf N}^{n-1}}
\sum_{m=1}^{\infty}ms_{K}s_{m}s_{m}^{*}s_{K}^{*}.\]
By definition, $s_{j}^{m}\Omega=(s_{j}^{*})^{m}\Omega=\Omega$
for any $m\geq 1$.
From these,
we obtain that $a_{n}a_{n}^{*}\Omega=j\Omega$ for any $n\in {\bf N}$.

It is sufficient to show ${\cal B}\Omega={\cal D}_{j}$.
By definition of the RBS, ${\cal B}\Omega\subset {\cal D}_{j}$.
We write $(a_{n}^{*})^{-1}\equiv a_{n}$
and $a_{n}^{0}=(a_{n}^{*})^{0}=I$ for convenience.
Then
for any $n\in {\bf N}$, there exists $M\in {\bf R}$ such that 
$s_{n}\Omega=M(a_{1}^{*})^{n-j}\Omega$.	
From this, we can derive that
\[s_{K}\Omega\in {\cal B}\Omega\quad (K\in {\bf N}^{*}).\]
Hence ${\cal D}_{j}\subset {\cal B}\Omega$.
Therefore the statement holds.

\noindent
(ii)
We see that $({\cal H},\pi_{12})$ is $P_{\infty}(12)$ with 
the GP vector $\Omega$.
The relations of $a_{n}$'s and $\Omega,\Omega^{'}$
are shown by assumption.
Let 
$V_{1}\equiv {\cal B}\Omega$ and $V_{2}\equiv {\cal B}\Omega^{'}$.
Then we see that 
$V_{1}$ and $V_{2}$ are $F_{12}$ and $F_{21}$, respectively.
By Lemma \ref{lem:bosonresult} (ii),
$V_{1}$ and $V_{2}$ are orthogonal in ${\cal H}$.

For $m\geq 1$ and $J=(j_{1},\ldots,j_{n})\in {\bf N}^{n}$,
\[
s_{J}\Omega=
\left\{\begin{array}{ll}
C_{n}a^{*(J-1)}a_{2}a_{4}\cdots a_{2m}\Omega\quad &(n=2m),\\
\\
C_{n}a^{*(J-1)}a_{1}a_{3}\cdots a_{2m-1}\Omega^{'}\quad & (n=2m-1)\\
\end{array}
\right.
\]
where 
$a^{*(J-1)}\equiv (a_{1}^{*})^{j_{1}-1}\cdots (a_{n}^{*})^{j_{n}-1}$ and 
$C_{n}\equiv \{(j_{1}-1)!\cdots (j_{n}-1)!\}^{-1/2}$.
From this,
$s_{J}\Omega\in V_{1}\oplus V_{2}$ for any $J\in {\bf N}^{*}$.
This implies that the reference subspace of ${\cal H}$
is a subspace of $V_{1}\oplus V_{2}$.
Hence $V_{1}\oplus V_{2}$ is dense in ${\cal H}$.

\noindent
(iii) 
From (i), (ii) and Lemma \ref{lem:bosonresult} (i),
we see that $\eta_{j}$ is $F_{j}$ $(j\geq 1)$,
$\eta_{12}$ is $F_{12}$ and $\eta_{21}$ is $F_{21}$.
From these and Lemma \ref{lem:bosonresult} (ii), the statement holds.

\noindent
(iv) 
From Lemma \ref{lem:bosonresult} (iii), the statement holds.
\qedh

%
%
\sftt{Example}
\label{section:third}

%
%
\ssft{Fock representation of RBS}
\label{subsection:thirdone}
From Theorem \ref{Thm:restriction} (i),
we obtain a correspondence between state vectors in the Bose-Fock space 
and vectors in the permutative representation $P_{\infty}(1)$ as follows:
%
%
\begin{equation}
\label{eqn:rbsfive}
(a^{*}_{1})^{j_{1}-1}\cdots (a^{*}_{k})^{j_{k}-1}\Omega
=\{(j_{1}-1)!\cdots (j_{k}-1)!\}^{1/2}s_{J}\Omega
\end{equation}
for $J=(j_{1},\ldots, j_{k})\in{\bf N}^{k}$.
This shows that any physical theory with the Bose-Fock space
is rewritten by $\coni$. Furthermore the Fock vacuum is interpreted
as the eigenvector of the generator $s_{1}$ of $\coni$.
For example, the one-particle state is given as follows:
\[a_{n}^{*}\Omega=s_{1}^{n-1}s_{2}\Omega\quad(n\geq 1).\]

On the other hand,
if the Fock representation of ${\cal B}$ is given,
then it is always extended to the action of $\coni$ as follows:
\[
\begin{array}{rl}
s_{m}\Omega =&\{(m-1)!\}^{-1/2}(a_{1}^{*})^{m-1}\Omega,\\
\\
s_{m}(a_{n_{1}}^{*})^{k_{1}}\cdots(a_{n_{p}}^{*})^{k_{p}}\Omega=&
\{(m-1)!\}^{-1/2}\,
(a_{1}^{*})^{m-1}(a_{n_{1}+1}^{*})^{k_{1}}\cdots(a_{n_{p}+1}^{*})^{k_{p}}\Omega,
\\
\\
s_{m}^{*}\Omega =&\delta_{m,1}\Omega,\\
\\
s_{m}^{*}(a_{n_{1}}^{*})^{k_{1}}\cdots(a_{n_{p}}^{*})^{k_{p}}\Omega
=&
\left\{
\begin{array}{ll}
\delta_{m,1}(a_{n_{1}-1}^{*})^{k_{1}}
\cdots(a_{n_{p}-1}^{*})^{k_{p}}\Omega\quad &(n_{1}\geq 2),\\
\\
\delta_{m,k_{1}+1}
\sqrt{k_{1}!}\,
(a_{n_{2}-1}^{*})^{k_{2}}\cdots(a_{n_{p}-1}^{*})^{k_{p}}\Omega\quad
& (n_{1}=1)\\
\end{array}
\right.
\end{array}
\]
for $1\leq n_{1}<\cdots<n_{p}$ and $k_{1},\ldots,k_{p}\in {\bf N}$.

%
%
\begin{ex}{\rm
Define the representation $(\ltn,\pi)$ of $\coni$ by
%
%
\begin{equation}
\label{eqn:exampleone}
\pi(s_{n})e_{m}\equiv e_{2^{n-1}(2m-1)}\quad(n,m\in {\bf N}).
\end{equation}
Then this is $P_{\infty}(1)$ with the GP vector $e_{1}$.
For the representation in (\ref{eqn:exampleone}),
the vacuum is $e_{1}$ and the subspace $H_{1}$ of one-particle states
is given by
\[H_{1}\equiv \overline{{\rm Lin}\langle\{e_{2^{n-1}+1}:n\geq 1\}\rangle}.\]
}
\end{ex}
We show that the above correspondence holds for $\con$
for any $2\leq N< \infty$.
%
%
\begin{prop}
\label{prop:restrictiontwo}
If we identify $\coni$ with a C$^{*}$-subalgebra of $\con$ 
by (\ref{eqn:embedding}),  then $P_{N}(1)|_{{\cal B}}=Fock$.
\end{prop}
%
%
\pr
Because $P_{N}(1)|_{\coni}=P_{\infty}(1)$,
$P_{N}(1)|_{{\cal B}}=P_{\infty}(1)|_{{\cal B}}=Fock$.
\qedh

Let $({\cal H},\pi)$ be $P_{N}(1)$ of $\con$ with the GP vector $\Omega$.
From (\ref{eqn:embedding}), 
the following holds for $1\leq n_{1}<n_{2}<\cdots<n_{m}$
and $k_{1},\ldots,k_{m}\in {\bf N}$:
%
%
\begin{equation}
\label{eqn:rbsnth}
(a_{n_{1}}^{*})^{k_{1}}\cdots (a_{n_{m}}^{*})^{k_{m}}\Omega
=\prod_{i=1}^{m}\sqrt{k_{i}!}\,\,
t_{1}^{n_{1}-1}t^{c_{1}-1}_{N}t_{b_{1}}T_{2}\cdots T_{m}\Omega
\end{equation}
where $T_{i}\equiv t_{1}^{n_{i}-n_{i-1}-1}t^{c_{i}-1}_{N}t_{b_{i}}$ 
for $i=2,\ldots,m$ and 
we define $c_{i}\in {\bf N}$ and $b_{i}\in\{1,\ldots,N-1\}$ by 
the equation $k_{i}=(N-1)(c_{i}-1)+b_{i}-1$.

%
%
\begin{ex}(Fock representation by $\co{2}$ and $\co{3}$)
\label{ex:otwo}
{\rm
From (\ref{eqn:rbsnth}), the following holds:
When $N=2$,
\[(a_{n_{1}}^{*})^{k_{1}}\cdots (a_{n_{m}}^{*})^{k_{m}}\Omega
=\prod_{i=1}^{m}\sqrt{k_{i}!}\,\,
t_{1}^{n_{1}-1}t^{k_{1}}_{2}t_{1}^{n_{2}-n_{1}}t^{k_{2}}_{2}
\cdots t_{1}^{n_{m}-n_{m-1}}t^{k_{m}}_{2}\Omega.\]
When $N=3$,
\[(a_{n_{1}}^{*})^{k_{1}}\cdots (a_{n_{m}}^{*})^{k_{m}}\Omega
=\prod_{i=1}^{m}\sqrt{k_{i}!}\,\,
t_{1}^{n_{1}-1}t^{c_{1}-1}_{3}t_{b_{1}}T_{2}\cdots T_{m}\Omega\]
where $T_{i}\equiv t_{1}^{n_{i}-n_{i-1}-1}t^{c_{i}-1}_{3}t_{b_{i}}$ 
for $i=2,\ldots,m$ and we define $c_{i}\in {\bf N}$ and $b_{i}\in\{1,2\}$ 
by $k_{i}=2(c_{i}-1)+b_{i}-1$.
}
\end{ex}

%
%
\ssft{Interpretation of representations by infinite product}
\label{subsection:thirdtwo}
In this subsection, we consider representations
$F_{j}$ $(j\geq 2)$, $F_{12}$ and $F_{21}$ of bosons
in Definition \ref{defi:bosonrep} from a viewpoint of Fock representation.
Formal infinite products of operators are introduced for this purpose.
%
%
\subsubsection{$F_{j}$}
\label{subsubsection:thirdtwoone}
For the cyclic vector $\Omega$ of $F_{j}$ in Definition \ref{defi:bosonrep}
with $j\geq 2$,
it seems that the formal vector 
%
%
\begin{equation}
\label{eqn:fj}
\Omega^{'}\equiv \left(\prod_{n=1}^{\infty}a_{n}^{j-1}\right)\,\Omega
\end{equation}
is a new vacuum of $F_{j}$ up to normalization constant.
The cyclic subspace by $\Omega^{'}$ is 
equivalent to the Fock representation because
$a_{n}\Omega^{'}=0$ for each $n$ by formal computation.
However such vector can not be defined in the representation
space of $F_{j}$.
Furthermore $F_{j}$ is not equivalent to the Fock representation $F_{1}$
when $j\ne 1$ by Lemma \ref{lem:bosonresult} (ii).
However, the formal notation (\ref{eqn:fj}) often appears
in theoretical physics and it excites curiosity.
If we regard that (\ref{eqn:fj}) is justified by $F_{j}$,
then (\ref{eqn:fj}) obtains a meaning of 
the operation in the representation theory. 

%
%
\subsubsection{$F_{12}$ and $F_{21}$}
\label{subsubsection:thirdtwotwo}
According to the case $F_{j}$, we write the Fock vacuum by 
the cyclic vector $\Omega$ of $F_{12}$. Then
we obtain the formal vector $\Omega^{'}$ as follows:
%
%
\begin{equation}
\label{eqn:onetwo}
\Omega^{'}\equiv \left(\prod_{n=1}^{\infty}a_{2n}\right)\,\Omega.
\end{equation}
Of course, $\Omega^{'}$ never be defined in the representation space $F_{12}$.  

In the same way,
we write the Fock vacuum by the cyclic vector $\Omega$ of $F_{21}$. Then
we obtain the formal vector $\Omega^{'}$ as follows:
%
%
\begin{equation}
\label{eqn:twoone}
\Omega^{'}\equiv \left(\prod_{n=1}^{\infty}a_{2n-1}\right)\,\Omega.
\end{equation}

\noindent
{\bf Acknowledgement:}
The author would like to Mitsuo Abe for 
his idea of the recursive boson system.

\appendix

\ww
{\Large {\bf Appendix}}

%
\sftt{Proof of Lemma \ref{lem:permutative}}
\label{section:appone}

\noindent
(i) 
Fix $j\geq 1$.
We introduce an orthonormal basis of 
a given representation belonging to $P_{\infty}(j)$.
Let $({\cal H},\pi)$ be $P_{\infty}(j)$ with the GP vector $\Omega$.
We simply denote $\pi(s_{n})$ by $s_{n}$ for each $n$.
Define the subset $\Lambda_{j}$ of ${\bf N}^{*}\equiv 
\coprod_{l\geq 1}{\bf N}^{l}$ by
\[\Lambda_{j}\equiv\{(m),J\cup (n),\,n,m\geq 1,\,n\ne j,\,J\in{\bf N}^{*}\}\]
and $v_{J}\equiv s_{J}\Omega$ for $J\in{\bf N}^{*}$.
Because $s_{j}\Omega=\Omega$, we see that
$\{s_{J}s_{K}^{*}\Omega:J,K\in{\bf N}^{*}\}=\{s_{J}\Omega:J\in\Lambda_{j}\}$.
Hence
${\rm Lin}\langle \{v_{J}:J\in\Lambda_{j}\}\rangle$ is dense in ${\cal H}$.
Furthermore $\langle v_{J}|\Omega\rangle=0$ when $J\ne (j)$.
This implies that 
$\langle v_{J}|v_{K}\rangle=\delta_{J,K}$ for $J,K\in\Lambda_{j}$.
In consequence
$\{v_{J}:J\in\Lambda_{j}\}$ is a complete orthonormal basis of ${\cal H}$.
The construction of $\{v_{J}:J\in\Lambda_{j}\}$
is independent of the choice of ${\cal H}$
except the existence of GP vector $\Omega$.
Hence $P_{\infty}(j)$ is uniquely up to unitary equivalence.

Assume that $({\cal H},\pi)$ is $P_{\infty}(12)$ with the GP vector of $\Omega$.
We identify $\pi(s_{n})$ with $s_{n}$ for each $n$.
By definition, we see that
$\{s_{J}\Omega:J\in {\bf N}^{*}\}$ spans a dense subspace of ${\cal H}$.
Define the  sequence $\{T_{n}\in {\bf N}^{*}:n\in {\bf N}\}$ of multiindices
by $T_{2k}\equiv (12)^{k}$ and $T_{2k-1}=(12)^{k-1}\cup (1)$ for each $k\geq 1$.
If $J\in {\bf N}^{n}$,  then
\[\langle v_{J}|\Omega\rangle =\delta_{J,T_{n}}.\]
From this,
the orthonormal basis $\{v_{J}:J\in\Lambda_{12}\}$ of ${\cal H}$ is given by
\[v_{J}\equiv s_{J}\Omega\quad(J\in\Lambda_{12})\]
where
$\Lambda_{12}\equiv \{(n2),(m),J\cup (k),J\cup (l2)
:n,m,k,l\in {\bf N},\,k\ne 2,\,l\ne 1,J\in{\bf N}^{*}\}$.
Hence the orthonormal basis of ${\cal H}$ is determined only by
the assumptions of $\Omega$.
Hence $P_{\infty}(12)$ is unique up to unitary equivalence.

\noindent
(ii)
Assume that $P_{\infty}(i)\sim P_{\infty}(j)$.
Then there exists a representation of $\coni$
with two cyclic vectors $\Omega$ and $\Omega^{'}$
satisfying $s_{i}\Omega=\Omega$ and $s_{j}\Omega^{'}=\Omega^{'}$.
Because $i\ne j$, $\langle\Omega|\Omega^{'}\rangle=0$.
Furthermore
we can verify that $\langle v_{J}|\Omega^{'}\rangle
=\delta_{J,(j)^{n}}\langle \Omega|\Omega^{'}\rangle=0$ 
for any $J\in\Lambda_{i}\cap {\bf N}^{n}$ 
with respect to the notation in the proof of (i) for $i$.
Hence $\langle v_{J}|\Omega^{'}\rangle=0$ for any $J\in\Lambda_{i}$.
This implies that 
$\Omega^{'}=0$. This contradicts with the choice of $\Omega^{'}$.
Therefore $P_{\infty}(i)\not\sim P_{\infty}(j)$.

Fix $i\geq 1$.
Assume that $P_{\infty}(12) \sim P_{\infty}(i)$.
Then there exists a representation of $\coni$
with two cyclic vectors $\Omega$ and $\Omega^{'}$ satisfying
$s_{12}\Omega=\Omega$ and $s_{i}\Omega^{'}=\Omega^{'}$.
Then $\langle \Omega|\Omega^{'}\rangle
=\langle s_{12}\Omega|s_{i}^{2}\Omega^{'}\rangle=0$.
For any $J\in {\bf N}^{n}$,
\[\langle v_{J}|\Omega^{'}\rangle
=\delta_{J,(i)^{n}}\langle \Omega|\Omega^{'}\rangle=0.\]
This implies $\Omega^{'}=0$.
This contradicts with the choice of $\Omega^{'}$.
Hence there exist no such cyclic vector. 
Therefore the statement holds.

\noindent
(iii)
We use the notation in the proof of (i).
Assume that $B\in {\cal B}({\cal H})$ satisfies
$[B,x]=0$ for any $x\in \coni$.
Then we can verify that
$\langle Bv_{J}|v_{K}\rangle=\delta_{J,K}\langle B\Omega|\Omega\rangle$
for each $J,K\in\Lambda_{j}$.
This implies that $B=\langle \Omega|B\Omega\rangle \cdot I\in {\bf C}I$.
Hence the statement holds.

Assume that $\coni$ acts on ${\cal H}$ and 
$\Omega \in {\cal H}$ is a cyclic vector such that $s_{12}\Omega=\Omega$.
Assume that $B\in {\cal B}({\cal H})$ satisfies
$[B,x]=0$ for any $x\in \coni$.
Then we can verify that
\[\langle Bv_{J}|v_{K}\Omega\rangle =\delta_{JK}\langle B\Omega|\Omega\rangle
\quad(J,K\in \Lambda_{12}).\]
From this, $B=\langle \Omega|B\Omega\rangle I\in {\bf C}I$.
Hence the statement holds.
\qedh

%
%
\sftt{Proof of Lemma \ref{lem:bosonresult}}
\label{section:apptwo}

\noindent
(i)
Fix $j\geq 1$.
By definition, the following is derived:
\[a_{n}^{j}\Omega=0,\quad
a_{n}^{k}(a_{n}^{*})^{k}\Omega=(j+k-1)\cdots j\Omega\quad (n,k\in {\bf N}).\]
In addition, if $j\geq 2$, then the following holds for $1\leq l\leq j-1$:
\[(a_{n}^{*})^{l}a_{n}^{l}\Omega=(j-1)\cdots (j-l)\Omega.\]
If $k\geq j$, then
$\langle \Omega|(a_{n}^{*})^{k}\Omega\rangle
=\langle a_{n}^{k}\Omega|\Omega\rangle=0$.
If $1\leq k\leq j-1$, then
\[\langle \Omega|(a_{n}^{*})^{k}\Omega\rangle
=C\langle \Omega|(a_{n}^{*})^{k}(a_{n}^{*})^{j-k}a_{n}^{j-k}\Omega\rangle
=C\langle a_{n}^{j}\Omega|a_{n}^{j-k}\Omega\rangle=0\]
where $C\equiv \{(j-1)\cdots k\}^{-1/2}$.
This implies 
$\langle \Omega|a_{n}^{k}\Omega\rangle=0$ when $1\leq k\leq j-1$.
In consequence, 
%
%
\begin{equation}
\label{eqn:orth}
\langle \Omega|a_{n}^{k}\Omega\rangle=
\langle \Omega|(a_{n}^{*})^{k}\Omega\rangle=0\quad(k,n\geq 1).
\end{equation}
From these, the family of the following vectors
is an orthonormal basis of the vector space ${\cal B}\Omega$:
%
%
\begin{equation}
\label{eqn:typej}
v=C\cdot (a_{n_{1}}^{*})^{k_{1}}\cdots (a_{n_{p}}^{*})^{k_{p}}
a_{m_{1}}^{l_{1}}\cdots a_{m_{q}}^{l_{q}}\Omega
\end{equation}
for $1\leq n_{1}<\cdots <n_{p}$ and $k_{1},\ldots,k_{p}\in {\bf N}$,
$1\leq m_{1}<\cdots <m_{q}$, $l_{1},\ldots, l_{q}\in \{1,\ldots,j-1\}$,
$\{n_{1},\ldots,n_{p}\}\cap \{m_{1},\ldots,m_{q}\}=\emptyset$
and $p,q\geq 0$ where we use notations $a_{n_{0}}^{*}=a_{m_{0}}=I$ and
\[C\equiv 
\left[\prod_{t=1}^{p}\{(j+k_{t}-1)\cdots j\}\cdot 
\prod_{r=1}^{q}\{(j-1)\cdots (j-l_{r})\}\right]^{-1/2}.\]
In particular, when $j=1$, we always assume $q=0$.
The existence of the canonical basis 
consisting of $v$'s in (\ref{eqn:typej}) implies
the uniqueness of the representation.
Therefore the statement holds for $F_{j}$. 

For the cyclic vector $\Omega$ of $F_{12}$, we see that
\[a_{2m-1}^{l}(a_{2m-1}^{*})^{l}\Omega=l!\,\Omega,\quad
a_{2m}^{l}(a_{2m}^{*})^{l}\Omega=(l+1)!\,\Omega,\quad a_{2m}^{2}\Omega=0\]
for $l,m\geq 1$.
Define
%
%
\begin{equation}
\label{eqn:onetwov}
v=C\cdot 
(a_{2n_{1}-1}^{*})^{k_{1}}\cdots 
(a_{2n_{p}-1}^{*})^{k_{p}}
(a_{2m_{1}}^{*})^{l_{1}}\cdots 
(a_{2m_{q}}^{*})^{l_{q}}
a_{2t_{1}}\cdots a_{2t_{r}}\Omega
\end{equation}
for
$1\leq n_{1}<\cdots<n_{p}$,
$1\leq m_{1}<\cdots <m_{q}$,
$1\leq t_{1}<\cdots <t_{r}$,
$\{m_{1},\ldots,m_{q}\}\cap \{t_{1},\ldots,t_{r}\}=\emptyset$
and $k_{1},\ldots,k_{p},l_{1},\ldots,l_{q}\in {\bf N}$
where
\[C=\{k_{1}!\cdots k_{p}!\, (l_{1}+1)!\cdots (l_{q}+1)!\}^{-1/2}.\]
Then the set of all such $v$'s in (\ref{eqn:onetwov})
is an orthonormal basis of ${\cal B}\Omega$.
Hence the uniqueness of $F_{12}$ holds.

For $F_{21}$, we can construct an orthonormal basis by replacing
the suffixes $2n$ and $2n-1$ in the proof of $F_{12}$.
Hence the uniqueness of $F_{21}$ holds.

\noindent
(ii)
Assume that $i\ne j$ and $F_{i}\sim F_{j}$.
Then there exists a representation of ${\cal B}$
with two cyclic vectors $\Omega$ and $\Omega^{'}$ satisfying 
%
%
\begin{equation}
\label{eqn:ineone}
a_{n}a_{n}^{*}\Omega=i\Omega,\quad
a_{n}a_{n}^{*}\Omega^{'}=j\Omega^{'}\quad(n\geq 1).
\end{equation}
From this, $\langle\Omega|\Omega^{'}\rangle=0$.
Let 
%
%
\begin{equation}
\label{eqn:vectorone}
x=(a_{n_{1}}^{*})^{k_{1}}\cdots (a_{n_{p}}^{*})^{k_{p}}
a_{m_{1}}^{l_{1}}\cdots a_{m_{q}}^{l_{q}}
\end{equation}
for $1\leq n_{1}<\cdots<n_{p}$, $1\leq m_{1}<\cdots<m_{q}$,
$k_{1},\ldots,k_{p}\in {\bf N}$, $l_{1},\ldots,l_{q}\in \{1,\ldots,j-1\}$
and $\{n_{1},\ldots,n_{p}\}\cap \{m_{1},\ldots,m_{q}\}=\emptyset$.
Define $M\equiv n_{p}+m_{q}+1$. 
Because $a_{M}a_{M}^{*}$ commutes $x$,
and $i\ne j$, $\langle x\Omega|\Omega^{'}\rangle=0$ from (\ref{eqn:ineone}).
From this and (\ref{eqn:typej}), $\Omega^{'}=0$.
This contradicts with the choice of $\Omega^{'}$.
Hence $F_{i}\not \sim F_{j}$ when $i\ne j$.

Assume that $F_{j}\sim F_{12}$ for some $j\geq 1$.
Then there exists a representation of ${\cal B}$
with two cyclic vectors $\Omega$ and $\Omega^{'}$ satisfying 
%
%
\begin{equation}
\label{eqn:inethree}
a_{n}a_{n}^{*}\Omega=j\Omega,\quad a_{2n-1}\Omega^{'}=0,\quad
a_{2n}^{*}a_{2n}\Omega^{'}=\Omega^{'}\quad(n\geq 1).
\end{equation}
From this, 
%
%
\begin{equation}
\label{eqn:inetwo}
a_{2n}a_{2n}^{*}\Omega^{'}=2\Omega^{'},\quad 
a_{2n-1}a_{2n-1}^{*}\Omega^{'}=\Omega^{'}\quad(n\geq 1).
\end{equation}
Hence $\langle\Omega|\Omega^{'}\rangle=0$.
Let $x$ be as in (\ref{eqn:vectorone}) and $M\equiv n_{p}+m_{q}+1$.
Because both $a_{2M-1}a_{2M-1}^{*}$ and $a_{2M}a_{2M}^{*}$ commute $x$,
$\langle x\Omega|\Omega^{'}\rangle=0$
from (\ref{eqn:inethree}) and (\ref{eqn:inetwo}).
This implies $\Omega^{'}=0$.
This contradicts with the choice of $\Omega^{'}$.
Therefore $F_{j}\not\sim F_{12}$ for any $j\geq 1$.
In the same way, 
we see that $F_{j}\not\sim F_{21}$ for any $j\geq 1$.

Assume $F_{12}\sim F_{21}$.
Then there exists a representation of ${\cal B}$ with
two cyclic vectors $\Omega$ and $\Omega^{'}$ satisfying
$a_{2n}^{*}a_{2n}\Omega=\Omega$ and $a_{2n-1}\Omega=0$ and 
$a_{2n-1}^{*}a_{2n-1}\Omega^{'}=\Omega^{'}$ and $a_{2n}\Omega^{'}=0$ 
for each $n\geq 1$. Then
$\langle \Omega|\Omega^{'}\rangle
=\langle a_{2n}^{*}a_{2n}\Omega|\Omega^{'}\rangle
=\langle a_{2n}^{*}\Omega|a_{2n}\Omega^{'}\rangle=0$.
Let 
%
%
\begin{equation}
\label{eqn:vectorthree}
x=(a_{2n_{1}-1}^{*})^{k_{1}}\cdots (a_{2n_{p}-1}^{*})^{k_{p}}
(a_{2m_{1}}^{*})^{l_{1}}\cdots (a_{2m_{q}}^{*})^{l_{q}}
a_{2t_{1}}\cdots a_{2t_{r}}
\end{equation}
and assume the assumption in (\ref{eqn:onetwov}) and $p+q+r\geq 1$.
Let $L\equiv 2n_{p}-1+2m_{q}+2t_{r}+1$.
Then 
\[\langle x\Omega|\Omega^{'}\rangle
=\langle x a_{2L}^{*}a_{2L}\Omega|\Omega^{'}\rangle
=\langle a_{2L}^{*}xa_{2L}\Omega|\Omega^{'}\rangle
=\langle xa_{2L}\Omega|a_{2L}\Omega^{'}\rangle=0.\]
This holds for any such $x$. 
Hence $\Omega^{'}=0$. 
This contradicts with the choice of
$\Omega^{'}$. Therefore Assume $F_{12}\not\sim F_{21}$.

\noindent
(iii)
Fix $j\geq 1$.
Let $\Omega$ be the cyclic vector $F_{j}$
such that $a_{n}a_{n}^{*}\Omega=j\Omega$ for each $n\in {\bf N}$.
Assume that $B\in {\cal B}({\cal H})$ satisfies 
$[a_{n},B]=[a_{n}^{*},B]=0$ for each $n$.
Let $x$ be as in (\ref{eqn:vectorone}) and assume $p+q\geq 1$.
Because $B$ commutes $x$,
$\langle B\Omega|x\Omega\rangle=0$ from (\ref{eqn:orth}).
This implies that
$\langle Bw|z\rangle=0$ for $w,z\in {\cal E},\,w\ne z$
where ${\cal E}$ is the family of $v$'s in (\ref{eqn:typej}).
Therefore the off-diagonal part of $B$
with respect to vectors in ${\cal E}$ is zero.
Furthermore, we obtain that
$\langle Bv|v\rangle=\langle B\Omega|\Omega\rangle$ for any $v\in {\cal E}$.
This implies that $B=\langle \Omega|B\Omega\rangle I\in {\bf C}I$.
Hence $F_{j}$ is irreducible.

Let ${\cal S}$ be the set of all vectors $v$ in (\ref{eqn:onetwov}).
We see that $\langle \Omega|Bv\rangle=0$ for $v\in {\cal S}\setminus\{\Omega\}$.
From this,  $\langle v|Bw \rangle=0$ for $v,w\in {\cal S}$, $v\ne w$.
Furthermore, from the form of $v\in {\cal S}$, we obtain that
$\langle v|Bv\rangle =\langle \Omega|B\Omega\rangle$ for any $v\in {\cal S}$.
Therefore $B=\langle \Omega|B\Omega\rangle I\in {\bf C}I$.
Hence $F_{12}$ is irreducible.
We can prove the irreducibility of $F_{21}$
by replacing the suffixes $2n$ and $2n-1$ in the proof of $F_{12}$.
Hence $F_{21}$ is also irreducible.
\qedh


\end{document}